\documentclass[11pt,twoside,a4paper]{article}
\usepackage{amssymb}
\usepackage{amsmath}
\usepackage[dvips]{graphicx}
\title{Rigidity and the Lower Bound Theorem for Doubly Cohen-Macaulay Complexes}
\author{Eran Nevo \footnote{Institute of Mathematics, The Hebrew
University, Jerusalem Israel, E-mail address:
eranevo@math.huji.ac.il}}
 \newtheorem{thm}{Theorem}[section]
 
  \newtheorem{cor}[thm]{Corollary}
 \newtheorem{lem}[thm]{Lemma}
\newtheorem{de}[thm]{Definition}
 
\newtheorem{prob}[thm]{Problem}
\newtheorem{conj}[thm]{Conjecture}

\begin{document}
\maketitle

\begin{abstract}
We prove that for $d\geq 3$, the $1$-skeleton of any
$(d-1)$-dimensional doubly Cohen-Macaulay (abbreviated $2$-CM)
complex is generically $d$-rigid. This implies that Barnette's
lower bound inequalities for boundary complexes of simplicial
polytopes (\cite{Barnette},\cite{Barnette0}) hold for every $2$-CM
complex of dimension $\geq 2$ (see Kalai \cite{Kalai-LBT}).
Moreover, the initial part $(g_0,g_1,g_2)$ of the $g$-vector of a
$2$-CM complex (of dimension $\geq 3$) is an $M$-sequence. It was
conjectured by Bj\"{o}rner and Swartz \cite{Swartz} that the
entire $g$-vector of a $2$-CM complex is an $M$-sequence.
\end{abstract}

\section{Introduction}
The $g$-theorem gives a complete characterization of the
$f$-vectors of boundary complexes of simplicial polytopes. It was
conjectured by McMullen in 1970 and proved by Billera and Lee
\cite{Billera-Lee} (sufficiency) and by Stanley \cite{St}
(necessity) in 1980. A major open problem in $f$-vector theory is
the $g$-conjecture, which asserts that this characterization holds
for all homology spheres. The open part of this conjecture is to
show that the $g$-vector of every homology sphere is an
$M$-sequence, i.e. it is the $f$-vector of some order ideal of
monomials. Based on the fact that homology spheres are doubly
Cohen-Macaulay (abbreviated $2$-CM) and that the $g$-vector of
some other classes of $2$-CM complexes is known to be an
$M$-sequence (e.g. \cite{Swartz}), Bj\"{o}rner and Swartz
\cite{Swartz} recently suspected that
\begin{conj}\label{g-conj 2CM}(\cite{Swartz}, a weakening of Problem 4.2.)
The $g$-vector of any $2$-CM complex is an $M$-sequence.
\end{conj}
We prove a first step in this direction, namely:
\begin{thm}\label{g_2-2CM}
Let $K$ be a $(d-1)$-dimensional $2$-CM simplicial complex (over
some field) where $d\geq 4$. Then $(g_0(K),g_1(K),g_2(K))$ is an
$M$-sequence.
\end{thm}
This theorem follows from the following theorem, combined with an
interpretation of rigidity in terms of the face ring
(Stanley-Reisner ring), due (implicitly) to Lee \cite{Lee}.
\begin{thm}\label{thmRigid-2CM}
Let $K$ be a $(d-1)$-dimensional $2$-CM simplicial complex (over
some field) where $d\geq 3$. Then $K$ has a generically $d$-rigid
$1$-skeleton.
\end{thm}
Kalai \cite{Kalai-LBT} showed that if a simplicial complex $K$ of
dimension $\geq 2$ satisfy the following conditions then it
satisfies Barnette's lower bound inequalities:

(a) $K$ has a generically $(dim(K)+1)$-rigid $1$-skeleton.

(b) For each face $F$ of $K$ of codimension $>2$, its link
$lk_K(F)$ has a generically $(dim(lk_K(F))+1)$-rigid $1$-skeleton.

(c) For each face $F$ of $K$ of codimension $2$, its link
$lk_K(F)$ (which is a graph) has at least as many edges as
vertices.

Kalai used this observation to prove that Barnette's inequalities
hold for a large class of simplicial complexes.

Observe that the link of a vertex in a $2$-CM simplicial complex
is $2$-CM, and that a $2$-CM graph is $2$-connected. Combining it
with Theorem \ref{thmRigid-2CM} and the above result of Kalai we
conclude:
\begin{cor}\label{LBT-2CM}
Let $K$ be a $(d-1)$-dimensional $2$-CM simplicial complex where
$d\geq 3$. For all $0\leq i\leq d-1$ $f_i(K)\geq f_i(n,d)$ where
$f_i(n,d)$ is the number of $i$-faces in a (equivalently every)
stacked $d$-polytope on $n$ vertices. (Explicitly,
$f_{d-1}(n,d)=(d-1)n-(d+1)(d-2)$ and
$f_{i}(n,d)=\binom{d}{i}n-\binom{d+1}{i+1}i$ for $1\leq i\leq
d-2$.) $\square$
\end{cor}
Theorem \ref{thmRigid-2CM} is proved by decomposing $K$ into a
union of minimal $(d-1)$-cycle complexes (Fogelsanger's notion
\cite{Fogelsanger}). Each of these pieces has a generically
$d$-rigid $1$-skeleton (\cite{Fogelsanger}), and the decomposition
is such that gluing the pieces together results in a complex with
a generically $d$-rigid $1$-skeleton. The decomposition is
detailed in Theorem \ref{structureLemma}.

This paper is organized as follows: In Section \ref{SecRigidity}
we give the necessary background from rigidity theory, explain the
connection between rigidity and the face ring, and reduce the
results mentioned in the Introduction to Theorem
\ref{structureLemma}. In Section \ref{SecDecomposition} we give
the necessary background on $2$-CM complexes, prove Theorem
\ref{structureLemma} and discuss related problems and results.
\\
\section{Rigidity}\label{SecRigidity}
The presentation of rigidity here is based mainly on the one in
Kalai \cite{Kalai-LBT}.
\newline
Let $G=(V,E)$ be a graph. A map $f:V\rightarrow \mathbb{R}^{d}$ is
called a $d-embedding$. It is $rigid$ if any small enough
perturbation of it which preserves the lengths of the edges is
induced by an isometry of $\mathbb{R}^{d}$. Formally, $f$ is
called $rigid$ if there exists an $\varepsilon>0$ such that if
$g:V\rightarrow \mathbb{R}^{d}$ satisfies
$d(f(v),g(v))<\varepsilon$ for every $v\in V$ and
$d(g(u),g(w))=d(f(u),f(w))$ for every $\{u,w\}\in E$, then
$d(g(u),g(w))=d(f(u),f(w))$ for every $u,w\in V$ (where $d(a,b)$
denotes the Euclidean distance between the points $a$ and $b$).

$G$ is called $generically\ d-rigid$ if the set of its rigid
$d$-embeddings is open and dense in the topological vector space
of all of its $d$-embeddings.

Let $V=[n]$, and let $Rig(G,f)$ be the $dn\times |E|$ matrix
which is defined as follows:
for its column corresponding to $\{v<u\}\in E$ put the vector
$f(v)-f(u)$ (resp. $f(u)-f(v)$) at the entries of the $d$ rows
corresponding to $v$ (resp. $u$) and zero otherwise. $G$ is generically
$d$-rigid iff $Im(Rig(G,f))=Im(Rig(K_{V},f)$ for a generic $f$,
where $K_{V}$ is the complete graph on $V$. $Rig(G,f)$ is called the $rigidity\ matrix$ of
$G$ (its rank is independent of the generic $f$ that we choose).

Let $G$ be the $1$-skeleton of a $(d-1)$-dimensional simplicial
complex $K$. We define $d$ generic degree-one elements in the
polynomial ring $A=\mathbb{R}[x_1,..,x_n]$ as follows:
$\Theta_i=\sum_{v\in [n]}f(v)_i x_v$ where $f(v)_i$ is the
projection of $f(v)$ on the $i$-th coordinate, $1\leq i\leq d$.
Then the sequence $\Theta=(\Theta_1,..,\Theta_d)$ is an l.s.o.p.
for the face ring $\mathbb{R}[K]=A/I_{K}$ ($I_K$ is the ideal in
$A$ generated by the monomials whose support is not an element of
$K$). Let $H(K)=\mathbb{R}[K]/(\Theta)=H(K)_0\oplus
H(K)_1\oplus...$ where $(\Theta)$ is the ideal in $A$ generated by
the elements of $\Theta$ and the grading is induced by the degree
grading in $A$. Consider the multiplication map $\omega:
H(K)_1\longrightarrow H(K)_2$, $m\rightarrow \omega m$ where
$\omega=\sum_{v\in [n]}x_v$. Lee \cite{Lee} proved that
\begin{equation}\label{eqLee}
dim_{\mathbb{R}} Ker(Rig(G,f)) = dim_{\mathbb{R}}
H(K)_2-dim_{\mathbb{R}} \omega(H(K)_1).
\end{equation}
Assume that $G$ is generically $d$-rigid. Then $dim_{\mathbb{R}}
Ker(Rig(G,f)) = f_1(K)-rank(Rig(K_{V},f)) = g_2(K) =
dim_{\mathbb{R}} H(K)_2-dim_{\mathbb{R}} H(K)_1$. Combining with
(\ref{eqLee}), the map $\omega$ is injective, and hence
$dim_{\mathbb{R}} (H(K)/(\omega))_i = g_i(K)$ for $i=2$; clearly
this holds for $i=0,1$ as well. Hence $(g_o(K),g_1(K),g_2(K))$ is
an $M$-sequence. We conclude that Theorem \ref{thmRigid-2CM}
implies Theorem \ref{g_2-2CM}, via the following algebraic result:
\begin{thm}\label{thm-w}
Let $K$ be a $(d-1)$-dimensional $2$-CM simplicial complex (over
some field) where $d\geq 3$. Then the multiplication map $\omega:
H(K)_1\longrightarrow H(K)_2$ is injective. $\square$
\end{thm}

In order to prove Theorem \ref{thmRigid-2CM}, we need the concept
of minimal cycle complexes, introduced by Fogelsanger
\cite{Fogelsanger}. We summarize his theory below.

Fix a field $k$ (or more generally, any abelian group) and
consider the formal chain complex on a ground set $[n]$,
$C=(\oplus\{kT :T\subseteq [n]\}, \partial)$, where
$\partial(1T)=\sum_{t\in T} sign(t,T)T\setminus \{t\}$ and
$sign(t,T)=(-1)^{|\{s\in T: s<t\}|}$. Define $subchain$, $minimal\
d-cycle$ and $minimal\ d-cycle\ complex$ as follows:
$c'=\sum\{b_TT: T\subseteq [n], |T|=d+1\}$ is a $subchain$ of a
$d$-chain $c=\sum\{a_TT: T\subseteq [n], |T|=d+1\}$ iff for every
such $T$, $b_T=a_T$ or $b_T=0$. A $d$-chain $c$ is a $d-cycle$ if
$\partial(c)=0$, and is a $minimal\ d-cycle$ if its only subchains
which are cycles are $c$ and $0$. A simplicial complex $K$ which
is spanned by the support of a $minimal\ d-cycle$ is called a
$minimal\ d-cycle\ complex$ (over $k$), i.e. $K=\{S: \exists T\
S\subseteq T, a_T\neq 0\}$ for some minimal $d$-cycle $c$ as
above. For example, triangulations of connected manifolds without
boundary are minimal cycle complexes - fix $k=\mathbb{Z}_2$ and
let the cycle be the sum of all facets.

The following is the main result in Fogelsanger's thesis.
\begin{thm}\label{thmFog}(Fogelsanger \cite{Fogelsanger})
For $d\geq 3$, every minimal $(d-1)$-cycle complex has a generically $d$-rigid $1$-skeleton.
\end{thm}

We will need the following gluing lemma, due of Asimov and Roth,
who introduced the concept of generic rigidity of graphs
\cite{Asi-Roth1}.
\begin{thm}\label{gluingLemma}(Asimov and Roth \cite{Asi-Roth2})
Let $G_1$ and $G_2$ be generically $d$-rigid graphs. If $G_1\cap
G_2$ contains at least $d$ vertices, then $G_1\cup G_2$ is
generically $d$-rigid.
\end{thm}
Now we are ready to conclude Theorem \ref{thmRigid-2CM} from the decomposition theorem, Theorem \ref{structureLemma}.
\newline
$proof\ of\ Theorem\ \ref{thmRigid-2CM}$: Consider a decomposition
sequence of $K$ as guaranteed by Theorem \ref{structureLemma},
$K=\cup^m_{i=1} S_i$. By Theorem \ref{thmFog} each $S_i$ has a
generically $d$-rigid $1$-skeleton. By Theorem \ref{gluingLemma}
for all $2\leq i\leq m$ $\cup^i_{j=1} S_j$ has a generically
$d$-rigid $1$-skeleton, in particular $K$ has a generically
$d$-rigid $1$-skeleton ($i=m$). $\square$
\newline
\textbf{Remark}: One can verify that Theorems \ref{thmFog} and
\ref{gluingLemma}, and hence also Theorem \ref{thmRigid-2CM},
continue to hold when replacing "generically $d$-rigid" by the
notion "$d$-hypperconnected", introduced by Kalai \cite{Kalai-56}.
Both of these assertions have an interpretation in terms of
algebraic shifting, introduced by Kalai (see e.g. his survey
\cite{skira}), namely: for both the exterior and symmetric
shifting operators over the field $\mathbb{R}$, denoted by
$\Delta$, $\{d,n\}\in \Delta(K)$. The existence of this edge in
the shifted complex implies the non-negativity of $g_2(K)$.

\section{Decomposing a $2$-CM complex}\label{SecDecomposition}
\begin{de}
A simplicial complex $K$ is $2-CM$ (over a fixed field $k$) if it
is Cohen-Macaulay and for every vertex $v\in K$, $K-v$ is
Cohen-Macaulay of the same dimension as $K$.
\end{de}
Here $K-v$ is the simplicial complex $\{T\in K: v\notin T\}$. By a
theorem of Reisner \cite{Reisner}, a simplicial complex $L$ is
Cohen-Macaulay iff it is pure and for every face $T\in L$
(including the empty set) and every $i<dim(lk_L(T)$,
$\tilde{H}_i(lk_L(T);k)=0$ where $lk_L(T)=\{S\in L: T\cap
S=\emptyset, T\cup S\in L\}$ and $\tilde{H}_i(M;k)$ is the reduced
$i$-th homology of $M$ over $k$. The proof of Theorem
\ref{structureLemma} is by induction on $dim(K)$. Let us first
consider the case where $K$ is $1$-dimensional.

A (simple finite) graph is $2$-\emph{connected} if after a
deletion of any vertex from it, the remaining graph is connected
and non trivial (i.e. is not a single vertex nor empty). Note that
a graph is $2$-CM iff it is $2$-connected.
\begin{lem}\label{graphLemma}
A graph $G$ is $2$-connected iff there exists a decomposition
$G=\cup^m_{i=1} C_i$ such that each $C_i$ is a simple cycle and
for every $1<i\leq m$, $C_i\cap(\cup_{j<i}C_j)$ contains an edge.

Moreover, for each $i_0\in [m]$ the $C_i$'s can be reordered by a
permutation $\sigma:[m]\rightarrow[m]$ such that
$\sigma^{-1}(1)=i_0$ and for every $i>1$,
$C_{\sigma^{-1}(i)}\cap(\cup_{j<i}C_{\sigma^{-1}(j)})$ contains an
edge.
\end{lem}
$Proof$: Whitney \cite{Whitney21} showed that a graph $G$ is
$2$-connected iff it has an open ear decomposition, i.e. there
exists a decomposition $G=\cup^m_{i=0} P_i$ such that each $P_i$
is a simple open path, $P_0$ is an edge, $P_0\cup P_1$ is a simple
cycle and for every $1<i\leq m$ $P_i\cap(\cup_{j<i}P_j)$ equals
the $2$ end vertices of $P_i$.

Assume that $G$ is $2$-connected and consider an open ear
decomposition as above. Let $C_1=P_0\cup P_1$. For $i>1$ choose a
simple path $\tilde{P}_i$ in $\cup_{j<i}P_j$ that connects the $2$
end vertices of $P_i$, and let $C_i=P_i\cup \tilde{P_i}$.
$(C_1,...,C_m)$ is the desired decomposition sequence of $G$.

Let $C$ be the graph whose vertices are the $C_i$'s and two of
them are neighbors iff they have an edge in common. Thus, $C$ is
connected, and hence the 'Moreover' part of the Lemma is proved.

The other implication, that such a decomposition implies
$2$-connectivity, will not be used in the sequel, and its proof is
omitted. $\square$
\\

For the induction step we need the following cone lemma. For $v$ a
vertex not in the support of a $(d-1)$-chain $c$, let $v*c$ denote
the following $d$-chain: if $c=\sum\{a_T T: v\notin T\subseteq
[n], |T|=d\}$ where $a_T\in k$ for all $T$, then
$v*c=\sum\{sign(v,T)a_T T\cup\{v\}: v\notin T\subseteq [n],
|T|=d\}$ where $sign(v,T)=(-1)^{|\{t\in T: t<v\}|}$.

\begin{lem}\label{sphereLemma}
Let $s$ be a minimal $(d-1)$-cycle and let $c$ be a minimal
$d$-chain such that $\partial(c)=s$, i.e. $c$ has no proper
subchain $c'$ such that $\partial(c')=s$. For $v$ a vertex not in
any face in $supp(c)$,the support of $c$, define
$\tilde{s}=c-v*s$. Then $\tilde{s}$ is a minimal $d$-cycle.
\end{lem}
$Proof$:
$\partial(\tilde{s})=\partial(c)-\partial(v*s)=s-(s-v*\partial(s))=0$
hence $\tilde{s}$ is a $d$-cycle. To show that it is minimal, let
$\hat{s}$ be a subchain of $\tilde{s}$ such that
$\partial(\hat{s})=0$. Note that $supp(c)\cap
supp(v*s)=\emptyset$.
\newline Case 1: $v$ is contained in a face in $supp(\hat{s})$. By the minimality
of $s$, $supp(v*s)\subseteq supp(\hat{s})$. Thus, by the
minimality of $c$ also $supp(c)\subseteq supp(\hat{s})$ and hence
$\hat{s}=\tilde{s}$.
\newline Case 2: $v$ is not contained in any face in $supp(\hat{s})$. Thus,
$supp(\hat{s})\subseteq supp(c)$. As $\partial(\hat{s})=0$ then
$\partial(c-\hat{s})=s$. The minimality of $c$ implies
$\hat{s}=0$. $\square$

\begin{thm}\label{structureLemma}
Let $K$ be a $d$-dimensional $2$-CM simplicial complex over a
field $k$ ($d\geq 1$). Then there exists a decomposition
$K=\cup^m_{i=1} S_i$ such that each $S_i$ is a minimal $d$-cycle
complex over $k$ and for every $i>1$, $S_i\cap(\cup_{j<i}S_j)$
contains a $d$-face.

Moreover, for each $i_0\in [m]$ the $S_i$'s can be reordered by a
permutation $\sigma:[m]\rightarrow[m]$ such that
$\sigma^{-1}(1)=i_0$ and for every $i>1$,
$S_{\sigma^{-1}(i)}\cap(\cup_{j<i}S_{\sigma^{-1}(j)})$ contains a
$d$-face.
\end{thm}

$proof$: The proof is by induction on $d$. For $d=1$, by Lemma
\ref{graphLemma} $K=\cup^{m(K)}_{i=1} C_i$ such that each $C_i$ is
a simple cycle and for every $i>1$ $C_i\cap(\cup_{j<i}C_j)$
contains an edge. Define $s_i=\sum\{sign_e(i)e: e\in (C_i)_1\}$,
then $s_i$ is a minimal $1$-cycle (orient the edges properly:
$sign_e(i)$ equals $1$ or $-1$ accordingly) whose support spans
the simplicial complex $C_i$. Moreover, by Lemma \ref{graphLemma}
each $C_{i_0}$, $i_0\in [m(K)]$, can be chosen to be the first in
such a decomposition sequence.

For $d>1$, note that the link of every vertex in a $2$-CM
simplicial complex is $2$-CM. For a vertex $v\in K$, as $lk_K(v)$
is $2$-CM then by the induction hypothesis
$lk_K(v)=\cup^{m(v)}_{i=1} C_i$ such that each $C_i$ is a minimal
$(d-1)$-cycle complex and for every $i>1$ $C_i\cap(\cup_{j<i}C_j)$
contains a $(d-1)$-face. Let $s_i$ be a minimal $(d-1)$-cycle
whose support spans $C_i$. As $K-v$ is CM of dimension $d$,
$\tilde{H}_{d-1}(K-v;k)=0$. Hence there exists a $d$-chain $c$
such that $\partial(c)=s_i$ and $supp(c)\subseteq K-v$.

Take $c_i$ to be such a chain with a support of minimal
cardinality. By Lemma \ref{sphereLemma}, $\tilde{s_i}=c_i-v*s_i$
is a minimal $d$-cycle. Let $S_i(v)$ by the simplicial complex
spanned by $supp(\tilde{s_i})$; it is a minimal $d$-cycle complex.
By the induction hypothesis, for every $i>1$
$S_i(v)\cap(\cup_{j<i}S_j(v))$ contains a $d$-face (containing
$v$). Thus, $K(v):=\cup_{j=1}^{m(v)}S_j(v)$ has the desired
decomposition for every $v\in K$. $K=\cup_{v\in Ver(K)}K(v)$ as
$st_K(v)\subseteq K(v)$ for every $v$, where $st_K(v)=\{T\in K:
T\cup\{v\}\in K\}$.

Let $v$ be any vertex of $K$. Since the $1$-skeleton of $K$ is
connected, we can order the vertices of $K$ such that $v_1=v$ and
for every $i>1$ $v_i$ is a neighbor of some $v_j$ where $1\leq
j<i$. Let $v_{l(i)}$ be such a neighbor of $v_i$. By the induction
hypothesis we can order the $S_j(v_i)$'s such that $S_1(v_i)$ will
contain $v_{l(i)}$, and hence, as $K$ is pure, will contain a
$d$-face which appears in $K(v_{l(i)})$ (this face contains the
edge $\{v_i,v_{l(i)}\}$). The resulting decomposition sequence
$(S_1(v_1),..,S_{m(v_1)}(v_1),S_1(v_2),..,S_{m(v_n)}(v_n))$ is as
desired.

Moreover, every $S_j(v_{i_0})$ where $i_0\in [n]$ and $j\in
[m(v_{i_0})]$ can be chosen to be the first in such a
decomposition sequence. Indeed, by the induction hypothesis
$S_j(v_{i_0})$ can be the first in the decomposition sequence of
$K(v_{i_0})$, and as mentioned before, the connectivity of the
$1$-skeleton of $K$ guarantees that each such prefix
$(S_1(v_{i_0}),..,S_{m(v_{i_0})}(v_{i_0}))$ can be completed to a
decomposition sequence of $K$ on the same $S_j(v_i)$'s. $\square$
\\

Theorem \ref{thmRigid-2CM} follows also from the following
corollary combined with Theorem \ref{thmFog}.
\begin{cor}\label{cor2CM->minCycle}
Let $K$ be a $d$-dimensional $2$-CM simplicial complex over a
field $k$ ($d\geq 1$). Then $K$ is a minimal cycle complex over
the Abelian group $\tilde{k}= k(x_1,x_2,...)$ whose elements are
finite linear combinations of the (variables) $x_i$'s with
coefficients in $k$.
\end{cor}
$Proof$: Consider a decomposition $K=\cup^m_{i=1} S_i$ as
guaranteed by Theorem \ref{structureLemma}, where
$S_i=\overline{supp(c_i)}$ (the closure w.r.t. inclusion of
$supp(c_i)$) for some minimal $d$-cycle $c_i$ over $k$. Define
$\tilde{c_i}=x_i c_i$, thus $\tilde{c_i}$ is a minimal cycle over
$\tilde{k}$. Define $\tilde{c}=\sum_{i=1}^m \tilde{c_i}$. Clearly
$\tilde{c}$ is a cycle over $\tilde{k}$ whose support spans $K$.
It remains to show that $\tilde{c}$ is minimal. Let $\tilde{c}'$
be a subchain of $\tilde{c}$ which is a cycle, $\tilde{c}'\neq
\tilde{c}$. We need to show that $\tilde{c}'=0$. Denote by
$\tilde{\alpha_T}$ ($\tilde{\alpha_T}'$) the coefficient of the
set $T$ in $\tilde{c}$ ($\tilde{c}'$) and by $\tilde{\alpha_T}(i)$
the coefficient of the set $T$ in $\tilde{c_i}$. If
$\tilde{\alpha_T}'=0$ then for every $i$ such that
$\tilde{\alpha_T}(i)\neq 0$, the minimality of $\tilde{c_i}$
implies that $\tilde{\alpha_F}'=0$ whenever
$\tilde{\alpha_F}(i)\neq 0$. By assumption, there exists a set
$T_0$ such that $\tilde{\alpha_{T_0}}'=0\neq
\tilde{\alpha_{T_0}}$. In particular, there exists an index $i_0$
such that $\tilde{\alpha_{T_0}}(i_0)\neq 0$, hence
$\tilde{\alpha_F}'=0$ whenever $\tilde{\alpha_F}(i_0)\neq 0$. As
$S_{i_0}\cap(\cup_{j<i_0}S_j)$ contains a $d$-face in case
$i_0>1$, repeated application of the above argument implies
$\tilde{\alpha_F}'=0$ whenever $\tilde{\alpha_F}(1)\neq 0$.
Repeated application of the fact that $S_{i}\cap(\cup_{j<i}S_j)$
contains a $d$-face for $i=2,3,..$ and of the above argument shows
that $\tilde{\alpha_F}'=0$ whenever $\tilde{\alpha_F}(i)\neq 0$
for some $1\leq i\leq m$, i.e. $\tilde{c}'=0$. $\square$
\newline

A pure simplicial complex has a \emph{nowhere zero flow} if there
is an assignment of integer non-zero wights to all of its facets
which forms a $\mathbb{Z}$-cycle. This generalizes the definition
of a nowhere zero flow for graphs (e.g. \cite{Seymour} for a
survey).
\begin{cor}\label{corZ2CM->nowhereZeroFlow}
Let $K$ be a $d$-dimensional $2$-CM simplicial complex over
$\mathbb{Q}$ ($d\geq 1$). Then $K$ has a nowhere zero flow.
\end{cor}
$Proof$: Consider a decomposition $K=\cup^m_{i=1} S_i$ as
guaranteed by Theorem \ref{structureLemma}. Multiplying by a
common denominator, we may assume that each
$S_i=\overline{supp(c_i)}$ for some minimal $d$-cycle $c_i$ over
$\mathbb{Z}$ (instead of just over $\mathbb{Q}$). Let $N$ be the
maximal $|\alpha|$ over all nonzero coefficients $\alpha$ of the
$c_i$'s, $1\leq i\leq m$. Let $\tilde{c}=\sum_{i=1}^m
(N^m)^{i}c_i$. $\tilde{c}$ is a nowhere zero flow for $K$; we omit
the details. $\square$


\begin{prob}\label{probSphereDecomp}
Can the $S_i$'s in Theorem \ref{structureLemma} be taken to be homology spheres?
\end{prob}
Yhonatan Iron and I proved (unpublished) the following lemma:
\begin{lem}
Let $K$, $L$ and $K\cap L$ be simplicial complexes of the same
dimension $d-1$. Assume that $K$ and $L$ are weak-Lefschetz, i.e.
that multiplication by a generic degree-one element $g$ in
$H=H(K),H(L)$, $g:H_{i-1}\longrightarrow H_i$, is injective  for
all $i\leq \lfloor{d/2}\rfloor$. If $K\cap L$ is CM then $K\cup L$
is weak-Lefschetz.
\end{lem}
In view of this lemma, if the intersections
$S_i\cap(\cup_{j<i}S_j)$ in Theorem \ref{structureLemma} can be
taken to be CM, and the $S_i$'s can be taken to be homology
spheres, then Conjecture \ref{g-conj 2CM} would be reduced to the
long standing $g$-conjecture for homology spheres. Can the
intersections be guaranteed to be CM?

\section*{Acknowledgments}
I would like to thank my advisor Gil Kalai, Anders Bj\"{o}rner and
Ed Swartz for helpful discussions. This research was done during
the author's stay at Institut Mittag-Leffler, supported by the ACE
network.

\end{document}